\documentclass[11pt]{amsart}
\usepackage{amsthm, amsfonts, amssymb, amsmath, url}
\usepackage{raleway}

\newcommand{\C}{\mathbb{C}}

\newcommand{\Q}{\mathbb{Q}}

\newtheorem{thm}{Theorem}
\newtheorem{lem}[thm]{Lemma}

\theoremstyle{definition}

\newtheorem*{rem}{Remark}

\newtheorem*{ack}{Acknowledgments}

\begin{document}
\title{On an identity of Sylvester}

\date{\today}
\author{Bogdan Nica}

\begin{abstract}  
We discuss an algebraic identity, due to Sylvester, as well as related algebraic identities and applications.
\end{abstract}

\address{\newline Department of Mathematical Sciences \newline Indiana University--Purdue University Indianapolis}
\email{bnica@iu.edu}

\maketitle

\section{Introduction}
This note is concerned with the following algebraic identity.

\begin{thm}[Sylvester's identity]\label{thm: syl}
Let $x_1,\dots,x_n$ be distinct elements of a field $\mathbb{F}$, where $n\geq 2$. Then, for any non-negative integer $d$, we have
\begin{align}\label{eq: syl}
\sum_{i=1}^n \frac{x^{d}_i}{\prod_{j: \: j\neq i\:} (x_i-x_j)}=\sum_{\lambda_1+\dots+\lambda_n=d-n+1} x_1^{\lambda_1}\dots x_n^{\lambda_n}.
\end{align} 
\end{thm}

Let us comment on the right-hand side of \eqref{eq: syl}. Firstly, $\lambda_1,\dots,\lambda_n$ are understood to be non-negative integers. Thus the right-hand side of \eqref{eq: syl} is a \emph{polynomial} function in $x_1,\dots,x_n$, whereas the left-hand side of \eqref{eq: syl} is, a priori, a \emph{rational} function in $x_1,\dots,x_n$. 

Secondly, we note that the right-hand side of \eqref{eq: syl} is an empty sum for $d<n-1$, respectively a single term of $1$ when $d=n-1$. This means that Sylvester's identity \eqref{eq: syl} includes, in low degrees, the following identity due to Euler:
\begin{align}\label{eq: eul}
\sum_{i=1}^n \frac{x^{d}_i}{\prod_{j: \: j\neq i\:} (x_i-x_j)}=\begin{cases}
0 & \textrm{ if } d=0,\dots,n-2, \\
1 & \textrm{ if } d=n-1.
\end{cases}
\end{align}

Thirdly, the right-hand side of \eqref{eq: syl} can be expressed concisely as a complete homogeneous symmetric functions, namely $h_{d-n+1}(x_1,\dots,x_n)$. Let us recall that the complete symmetric function of degree $k$ in $n$ variables is given by
\[h_{k}(x_1,\dots,x_n)=\sum_{\lambda_1+\dots+\lambda_n=k} x_1^{\lambda_1}\dots x_n^{\lambda_n},\]
where, by convention, $h_{0}(x_1,\dots,x_n)=1$, and $h_{k}(x_1,\dots,x_n)=0$ for $k<0$.

The attribution of the identity \eqref{eq: syl} to Sylvester is explained by Bhatnagar in \cite[Sec.1]{Bhat}. It seems however that the identity \eqref{eq: syl}, and the Sylvester attribution, are not quite absorbed in the mathematical canon. On the one hand, Knuth \cite[p.472-473, Notes to Exer.33]{Knu} states both the identity and the attribution; Stanley \cite[Exer.7.4 on p.450]{S} gives the identity \eqref{eq: syl}, but attributes it to Jacobi--a precursor of Sylvester \cite[p.490]{S}. On the other hand, Sylvester's identity was rediscovered by Louck, cf. \cite[App.A]{LB} and \cite[Thm.2.2]{CL}, and then again by Strehl and Wilf \cite[Sec.4]{SW}. 

As a further illustration of this point, that Sylvester's identity is relatively unknown, consider the following result from a recent and very interesting work of Garcia, Omar, O'Neill, and Yih \cite[Thm.3]{G++}: for distinct $x_1\dots,x_n\in \C$, and $z\in \C$, it holds that
\begin{align*}
\sum_{p=0}^\infty \frac{h_p(x_1,\dots,x_n)}{(p+n-1)!}\: z^{p+n-1}=\sum_{i=1}^n \frac{e^{x_i z}}{\prod_{j: \: j\neq i\:} (x_i-x_j)}.
\end{align*} 
This is a key identity for the purposes of \cite{G++}. The proof given therein is quite lengthy, occupying well over three pages, and it ultimately relies on a more sophisticated apparatus involving Schur functions. The following three-line manipulation shows that the above identity is, essentially, a restatement of \eqref{eq: syl} as a power series identity:
\begin{align*}
\sum_{i=1}^n \frac{e^{x_i z}}{\prod_{j: \: j\neq i\:} (x_i-x_j)}&=\sum_{i=1}^n \frac{1}{\prod_{j: \: j\neq i\:} (x_i-x_j)} \bigg(\sum_{d=0}^\infty \frac{(x_iz)^d}{d!}\bigg)\\
&=\sum_{d=0}^\infty \bigg(\sum_{i=1}^n \frac{x_i^d}{\prod_{j: \: j\neq i\:} (x_i-x_j)}\bigg)\: \frac{z^d}{d!}\\
&=\sum_{d=n-1}^\infty h_{d-n+1}(x_1,\dots,x_n)\: \frac{z^d}{d!}.
\end{align*} 

Our own interest in Sylvester's identity was prompted by a recent work \cite{N}. Specifically, Sylvester's identity is used as a key ingredient in the proof of the so-called Complete Coefficient Theorem \cite[Thm.6.2]{N}. The present note expands an appendix to the preprint version of \cite{N} (\url{arXiv:2110.05616}). For other occurrences of Sylvester's identity, the reader is invited to follow the trail of references indicated by Bhatnagar \cite[Sec.3]{Bhat} and by Milne \cite[p.7]{M}. See also Section~\ref{sec: dill} herein for a concrete application of Sylvester's identity. 

This note is written with two goals in mind. We think that Sylvester's identity deserves to be better known, and the first goal is to promote it. The second goal is to present three different arguments for Sylvester's identity--all of them seemingly new. The first argument, given in Section~\ref{sec: 1}, is probably what Euler would have written had he considered the higher degree case of \eqref{eq: eul}. The second and the third argument are given in Section~\ref{sec: 2+3}. In Section~\ref{sec: exteul}, we extend Euler's identity \eqref{eq: eul} in a different direction.

Here is a quick repertoire of other approaches to Sylvester's identity that can be found in the literature. Both Strehl - Wilf \cite{SW} and Chen - Louck \cite{CL} sketch elegant proofs of \eqref{eq: syl} by means of generating functions. Another proof can be found in \cite[Sec.2]{Bhat}. Stanley \cite[p.490]{S} suggests a proof of \eqref{eq: syl} via a more general identity pertaining to Schur functions. Over the complex numbers, Sylvester's identity can also be proved using the residue theorem; see Louck - Biedenharn \cite{LB} and Knuth \cite{Knu}. 

\begin{ack}
I thank Darij Grinberg for thoughtful comments on a preliminary version of this note.
\end{ack}

\section{An application: Dilcher's identity}\label{sec: dill}
The nice identity
\begin{align*}
\sum_{i=1}^n (-1)^{i-1} \binom{n}{i} \frac{1}{i}=\sum_{i=1}^n \frac{1}{i}
\end{align*}
goes back to Euler. Dilcher \cite[Cor.3]{D} obtained the following generalization: for any positive integer $d$, it holds that
\begin{align}\label{eq: dill}
\sum_{i=1}^n (-1)^{i-1} \binom{n}{i} \frac{1}{i^d}=\sum_{1\leq i_1\leq i_2\leq \dots \leq i_d\leq n}\: \frac{1}{i_1i_2\dots i_d}.
\end{align}
Dilcher arrives at \eqref{eq: dill} indirectly, as the limiting case of certain $q$-identities. Let us show how  \eqref{eq: dill} can be obtained directly from Sylvester's identity. This approach appears to be new.

Working in the rational field $\Q$, we use \eqref{eq: syl} on $x_i=1/i$, for $i=1,\dots,n$, and with $d+n-1$ in place of $d$. A quick computation reveals that
\[\frac{x^{d+n-1}_i}{\prod_{j: \: j\neq i\:} (x_i-x_j)}=(-1)^{i-1} \binom{n}{i} \frac{1}{i^d}\]
for each $i$. So the left-hand side of \eqref{eq: syl} turns, in the particular case at hand, into the left-hand side of \eqref{eq: dill}. The right-hand side of \eqref{eq: syl} is
\[
h_d\left(1,\frac{1}{2},\dots,\frac{1}{n}\right)=\sum_{\lambda_1+\dots+\lambda_n=d}\: \frac{1}{1^{\lambda_1} \dots n^{\lambda_n}},
\]
and the latter sum is a rewriting of the right-hand side of \eqref{eq: dill}.

\section{Euler's identity}\label{sec: eul}
The identity \eqref{eq: eul} is first documented in a letter of Euler to Goldbach, dated September 25th 1762 \cite[p.1123-1124]{E}. Euler goes on to say that the identity appears ``to be more than a little curious; however, it seems to me that you may have been so kind as to communicate something similar to me a long time ago.'' So it may be that Goldbach is the true originator. In a subsequent letter, dated November 9th 1762, Euler proves the identity \eqref{eq: eul} by using a partial fraction decomposition. 

A modern approach to Euler's identity is via Lagrange interpolation. Recall, this says the following: given $n$ distinct inputs $x_1,\dots,x_n\in \mathbb{F}$, and $n$ outputs $y_1,\dots,y_n\in \mathbb{F}$, there exists a unique polynomial $P\in \mathbb{F}[X]$ of degree less than $n$ satisfying $P(x_i)=y_i$ for $i=1,\dots,n$. The following version is more pertinent for our needs.

\begin{lem}[Lagrange interpolation]
Let $x_1,\dots,x_n$ be distinct elements of a field $\mathbb{F}$. Then every polynomial $P\in \mathbb{F}[X]$ of degree less than $n$ satisfies the polynomial identity
\begin{align}\label{eq: LagInt}
P(X)=\sum_{i=1}^n P(x_i) \prod_{ j : \: j\neq i\:} \frac{X-x_j}{x_i-x_j}.
\end{align}
\end{lem}

In particular, for each $d=0,\dots, n-1$ we have
\begin{align}\label{eq: LagInt2}
X^d=\sum_{i=1}^n x_i^d \prod_{ j : \: j\neq i\:} \frac{X-x_j}{x_i-x_j},
\end{align}
Now the weighted power sum appearing in Euler's identity \eqref{eq: eul} is precisely the coefficient of $X^{n-1}$ on the right-hand side of \eqref{eq: LagInt2}; but, visibly, the coefficient of $X^{n-1}$ on the left-hand side of \eqref{eq: LagInt2} is $0$ when $d=0,\dots,n-2$, respectively $1$ when $d=n-1$. This completes the proof of \eqref{eq: eul}.

\section{First proof of Sylvester's identity}\label{sec: 1} 
We aim to obtain a recurrence relation for the rational function
\begin{align}\label{eq: wps}
S_d(x_1,\dots,x_n)=\sum_{i=1}^n \frac{x^{d}_i}{\prod_{j: \: j\neq i\:} (x_i-x_j)}.
\end{align} 
We start, very much like Euler did for his proof of the identity \eqref{eq: eul}, with the partial fraction decomposition 
\begin{align}\label{eq: epfd}
\frac{1}{(X-x_1)\dots(X-x_{n-1})}=\sum_{i=1}^{n-1}\Bigg( \prod_{\substack{j=1\\j\neq i}}^{n-1}\frac{1}{x_i-x_j}\Bigg) \frac{1}{X-x_i}.
\end{align}
Evaluating \eqref{eq: epfd} at the remaining node, $x_n$, we infer that
\begin{align}\label{eq: pfd}
 \prod_{\substack{j=1\\j\neq n}}^{n}\frac{1}{x_n-x_j}=-\sum_{i=1}^{n-1} \prod_{\substack{j=1\\j\neq i}}^{n}\frac{1}{x_i-x_j}.
\end{align}
Multiplying \eqref{eq: pfd} through by $x^d_n$ yields
\begin{align*}
\Bigg( \prod_{\substack{j=1\\j\neq n}}^{n}\frac{1}{x_n-x_j}\Bigg)x_n^d=-\sum_{i=1}^{n-1}\Bigg( \prod_{\substack{j=1\\j\neq i}}^{n}\frac{1}{x_i-x_j}\Bigg) x_n^d
\end{align*}
and so, adding $n-1$ terms on both sides, 
\begin{align*}
\sum_{i=1}^{n}\Bigg( \prod_{\substack{j=1\\j\neq i}}^{n}\frac{1}{x_i-x_j}\Bigg) x_i^d=\sum_{i=1}^{n-1}\Bigg( \prod_{\substack{j=1\\j\neq i}}^{n}\frac{1}{x_i-x_j}\Bigg) (x_i^d-x_n^d).
\end{align*}
The left-hand side of the above identity is $S_d(x_1,\dots,x_n)$. The right-hand side can be rewritten, after dividing $x_i^d-x_n^d$ by $x_i-x_n$, as 
\begin{align*}
\sum_{i=1}^{n-1}\Bigg( \prod_{\substack{j=1\\j\neq i}}^{n-1}\frac{1}{x_i-x_j}\Bigg) \sum_{e=0}^{d-1} x_i^e x_n^{d-e-1}=\sum_{e=0}^{d-1}x_n^{d-e-1}\sum_{i=1}^{n-1}\Bigg( \prod_{\substack{j=1\\j\neq i}}^{n-1}\frac{1}{x_i-x_j}\Bigg)  x_i^e. 
\end{align*}
We note that, on the right-hand side, the inner sum is $S_e(x_1,\dots,x_{n-1})$. All in all, we have obtained the recurrence 
\begin{align}\label{eq: recpfd}
S_d(x_1,\dots, x_n)=\sum_{e=0}^{d-1} S_e(x_1,\dots,x_{n-1}) \: x_{n}^{d-e-1}.
\end{align}
On the other hand, the complete symmetric function $h_{d-n+1}(x_1,\dots,x_n)$ satisfies the same recurrence relation. Indeed:
\begin{align*}
h_{d-n+1}(x_1,\dots,x_n)&=\sum_{\substack{\lambda_1+\dots+\lambda_n=d-n+1}} x_1^{\lambda_1}\dots x_n^{\lambda_n}\\
&=\sum_{e=n-2}^{d-1}\bigg(\sum_{\substack{\lambda_1+\dots+\lambda_{n-1}=e-n+2}} x_1^{\lambda_1}\dots x_{n-1}^{\lambda_{n-1}}\bigg) x_n^{d-e-1}\\
&=\sum_{e=0}^{d-1} h_{e-n+2}(x_1,\dots,x_{n-1}) \: x_{n}^{d-e-1}.
\end{align*}
The last step exploits the fact that $h_{e-n+2}(x_1,\dots,x_{n-1})=0$ for $e<n-2$.

The desired identity, $S_d(x_1,\dots,x_n)=h_{d-n+1}(x_1,\dots,x_n)$ for all integers $d\geq 0$ and $n\geq 2$, easily follows from these recurrences. In fact, we are afforded the choice of inducting on $d$, or on $n$. Let us check the base cases. The base case $d=0$ amounts to $S_0(x_1,\dots,x_n)=0$; this is, in effect, a rewriting of \eqref{eq: pfd}. The base case $n=2$, amounts to
\[\frac{x^{d}_1}{x_1-x_2}+\frac{x^{d}_2}{x_2-x_1}=\sum_{\lambda_1+\lambda_2=d-1} x_1^{\lambda_1} x_2^{\lambda_2},\]
which is obviously true. 

\begin{rem}
The above proof works for any non-negative integer $d$. In particular, it also accounts for Euler's identity \eqref{eq: eul}. Although Sylvester did not record a proof of the identity \eqref{eq: syl}, previous writings suggest that he might have also started from the partial fraction decomposition \eqref{eq: epfd}, cf. \cite[p.433]{Bhat}. 

The proof is historically motivated and circumvents Lagrange interpolation. It requires, in exchange, a related algebraic fact--the partial fraction decomposition \eqref{eq: epfd}. In fact, by applying \eqref{eq: LagInt} to the constant polynomial $1$, we exactly get \eqref{eq: epfd}.
\end{rem}

\section{Second and third proof of Sylvester's identity}\label{sec: 2+3} 
In the next two proofs, we take Euler's identity \eqref{eq: eul} for granted and we only deal with the higher degree case, $d\geq n$. Both proofs involve the elementary symmetric functions, given by
\[e_k(x_1,\dots,x_n)=\sum_{1\leq i_1<\dots<i_k\leq n} x_{i_1}\dots x_{i_k}\]
for $1\leq k\leq n$. By convention, $e_0(x_1,\dots,x_n)=1$, and $e_k(x_1,\dots,x_n)=0$ for $k>n$.

Our arguments exploit two key facts pertaining to the elementary symmetric functions. The first one is the fundamental relation between the two families of symmetric homogeneous functions, the complete and the elementary:
\begin{align}\label{eq: 2nd2}
h_d=\sum_{k=1}^n (-1)^{k-1} e_k h_{d-k}
\end{align}
for any positive integer $d$; see, for instance, \cite[Eq.(7.13) on p.296]{S}. In \eqref{eq: 2nd2}, and for the rest of this section, we lighten the notation by suppressing the variables $x_1,\dots,x_n$. 

The second one is the Vieta expansion formula:
\begin{align}\label{eq: vie}
\prod_{i=1}^n \:(X-x_i)=\sum_{k=0}^n\: (-1)^k e_k X^{n-k}.
\end{align}

\begin{proof}[Second proof of Sylvester's identity] We note the following general fact: weighted power sums of the form 
\[
W_d=\sum_{i=1}^n w_i\: x^d_i, 
\]
where each $w_i$ is a rational expression in $x_1,\dots,x_n$, satisfy a recurrence relation of order $n$
\begin{align}\label{eq: 2nd1}
W_{d}=\sum_{k=1}^n (-1)^{k-1} e_k W_{d-k}
\end{align}
for each $d\geq n$.

Indeed, by evaluating the Vieta expansion \eqref{eq: vie} at $x_i$, we deduce that
\[x_i^n=\sum_{k=1}^n\: (-1)^{k-1} e_k x_i^{n-k}\] for each $i$. We multiply each such relation by $w_i x_i^{d-n}$, where $d\geq n$, and we add up. The outcome is the recurrence \eqref{eq: 2nd1}.

Consider now the weighted power sums $S_d$ defined by \eqref{eq: wps}; we already know that $S_d=h_{d-n+1}$ for $d=0,\dots, n-1$ thanks to ~\eqref{eq: eul}. The recurrences \eqref{eq: 2nd1} and \eqref{eq: 2nd2} imply, by induction on $d$, that $S_d=h_{d-n+1}$ for each $d\geq n$. \end{proof}

\begin{proof}[Third proof of Sylvester's identity]
Let $d\geq n$, and consider the polynomial 
\[Q(X)=(X^{d-n}+h_1 X^{d-n-1}+\dots+h_{d-n})\prod_{i=1}^n \:(X-x_i).\]
Using, once again, the Vieta expansion for the latter product, we see that
\begin{align*}
Q(X)&=\Big(\sum_{j=0}^{d-n} h_j X^{d-n-j}\Big)\Big(\sum_{k=0}^n (-1)^k e_k X^{n-k}\Big)\\
&=\sum_{\ell=0}^d \Big(\sum_{\substack{0\leq k\leq n\\ 0\leq j\leq d-n\\ k+j=\ell}} (-1)^k e_k h_j \Big)X^{d-\ell}.
\end{align*}
The coefficient of $X^{d-\ell}$ is $1$ when $\ell=0$; in the range $0<\ell\leq d-n$, it is the full sum $\sum_{k=0}^n (-1)^k e_k h_{\ell-k}$, which vanishes according to \eqref{eq: 2nd2}. In the remaining range $d-n<\ell\leq d$, we reindex by setting $\ell=d-m$ for $m=0,\dots,n-1$. We conclude that $Q(X)=X^d-R(X)$, where
\begin{align}\label{eq: rem}
R(X)=\sum_{m=0}^{n-1} \Big(\sum_{k=n-m}^{n} (-1)^{k-1} e_k h_{d-m-k} \Big)X^{m}.
\end{align}
In the above formula, we allow for the possibility that $d-m-k<0$, for in that case $h_{d-m-k}=0$. Note that the degree of $R(X)$ is less than $n$, and that the coefficient of $X^{n-1}$ is 
\[\sum_{k=1}^{n} (-1)^{k-1} e_k h_{d-n+1-k}=h_{d-n+1}.\]

By writing $X^d=Q(X)+R(X)$, we see that the explicit polynomial $R(X)$ can be nicely interpreted as the remainder obtained upon dividing $X^d$ by $\prod_{i=1}^n \:(X-x_i)$. Equivalently, $R(X)$ is the interpolating polynomial of degree less than $n$ satisfying $R(x_i)=x_i^d$ for $1\leq i\leq n$. Applying now \eqref{eq: LagInt} to the polynomial $R(X)$, we obtain
\begin{align}\label{eq: higher eul}
R(X)=\sum_{i=1}^n x_i^d \prod_{ j : \: j\neq i\:} \frac{X-x_j}{x_i-x_j}.
\end{align}

The desired identity \eqref{eq: syl} follows by reading off the coefficient of $X^{n-1}$ on both sides. \end{proof}

\section{An extension of Euler's identity}\label{sec: exteul}
Let us consider again the interpolation identity \eqref{eq: LagInt2}. We have obtained Euler's identity \eqref{eq: eul} by reading off the leading coefficients; what identities are encoded in the remaining coefficients? According to Vieta's formula, for each $i=1,\dots, n$ we may expand 
\[\prod_{ j : \: j\neq i\:} (X-x_j)=\sum_{m=0}^{n-1}(-1)^m e_m(x_1,\dots, \widehat{x}_i, \dots,x_n)\: X^{n-1-m},\]
where the hat denotes omission. Consequently, the coefficient of $X^{n-1-m}$ on the right-hand side of \eqref{eq: LagInt2} equals 
\[(-1)^m \sum_{i=1}^n \frac{x_i^d\: e_m(x_1,\dots, \widehat{x}_i, \dots,x_n)}{\prod_{ j\neq i\:} (x_i-x_j)}.\]
On the left-hand side of \eqref{eq: LagInt2}, the coefficient of $X^{n-1-m}$ is $1$ if $n-1-m=d$, and $0$ otherwise. We deduce the following \emph{extended Euler identity}.

\begin{thm}
Let $x_1,\dots,x_n$ be distinct elements of a field $\mathbb{F}$, where $n\geq 2$. For each $0\leq d\leq n-1$ and $0\leq m\leq n-1$, we have
\begin{align}\label{eq: ME}
\sum_{i=1}^n \frac{x_i^d\: e_m(x_1,\dots, \widehat{x}_i, \dots,x_n)}{\prod_{ j: \: j\neq i\: } (x_i-x_j)} =\begin{cases}
(-1)^m & \textrm{ if } d+m=n-1,\\
0 & \textrm{ otherwise}.
\end{cases}
\end{align}
\end{thm}

We note that the usual Euler identity \eqref{eq: eul} is the case $m=0$ of \eqref{eq: ME}.

As an application of the extended Euler identity \eqref{eq: ME}, we prove that
\begin{align}\label{eq: F2}
\sum_{i=1}^n \prod_{\: j : \: j\neq i\:} \frac{1-a x_ix_j}{x_i-x_j}=\begin{cases}
a^{(n-1)/2} & \textrm{ if } n \textrm{ is odd},\\
0 & \textrm{ if } n \textrm{ is even}.
\end{cases}
\end{align}
for any parameter $a\in \mathbb{F}$. The case $a=1$ is a problem shortlisted for the 2019 edition of the International Mathematical Olympiad \cite{IMO}. The case $a=0$ is covered by Euler's identity \eqref{eq: eul}.

The expansion
\[\prod_{ j : \: j\neq i\:} (1-ax_ix_j)=\sum_{m=0}^{n-1}(-ax_i)^m\: e_m(x_1,\dots, \widehat{x}_i, \dots,x_n)\]
allows us to rewrite the left-hand side of \eqref{eq: F2} as
\[\sum_{i=1}^n \prod_{\: j : \: j\neq i\:} \frac{1-ax_ix_j}{x_i-x_j}=\sum_{m=0}^{n-1} \sum_{i=1}^n \frac{(-ax_i)^{m}\: e_m(x_1,\dots, \widehat{x}_i, \dots,x_n)}{\prod_{ j: \: j\neq i\: } (x_i-x_j)}.
\]
By \eqref{eq: ME}, the inner sum vanishes except when $2m=n-1$. Thus, the above double sum vanishes when $n$ is even; when $n$ is odd, it has a single term, corresponding to $m=(n-1)/2$, and that term equals $(-a)^m(-1)^m=a^m$. This completes the proof of \eqref{eq: F2}.

A higher degree version of the extended Euler identity \eqref{eq: ME}--in other words, an extended Sylvester's identity--can be extracted without difficulty from ~\eqref{eq: higher eul}. We leave that to the keen reader.


\end{document}